\documentclass{article}

\usepackage{myarxiv}
\usepackage[utf8]{inputenc} % allow utf-8 input
\usepackage[T1]{fontenc}    % use 8-bit T1 fonts
\usepackage{hyperref}       % hyperlinks
\usepackage{url}            % simple URL typesetting
\usepackage{booktabs}       % professional-quality tables
\usepackage{amsfonts}       % blackboard math symbols
\usepackage{nicefrac}       % compact symbols for 1/2, etc.
\usepackage{microtype}      % microtypography
\usepackage{lipsum}

\usepackage{amssymb,amsmath,amsfonts,amsthm,enumerate}
\usepackage{subfigure}
\usepackage{graphicx}
\usepackage{epsfig}
\usepackage{float}
\usepackage[usenames]{color}
\usepackage{cite}
\usepackage{bm}
\usepackage{enumitem}
\setlist{nosep}

\newtheorem{coro}{Corollary}

\newtheorem*{prob}{Inverse Problem}

\newtheorem{thm}{Theorem}

\allowdisplaybreaks
\def\f{\frac}

\def\R{\mathbb{R}}
\def\N{\mathbb{N}}
\def\C{\mathbb{C}}

\def\cU{\mathcal{U}}

\def\al{\alpha}
\def\be{\beta}

\def\ep{\epsilon}

\def\la{\lambda}
\def\La{\Lambda}

\def\na{\nabla}

\def\Om{\Omega}

\def\Sg{\Sigma}

\def\OOO{\Omega}
\def\ooo{\overline}
\def\ppp{\partial}
\def\ov{\overline}
\def\pa{\partial}

\def\wh{\widehat}

\def\pppa{\partial_t^{\alpha}}

\def\va{\varphi}

\title{Inverse Problems of Determining Coefficients of the Fractional Partial Differential Equations}

\author{
  Zhiyuan Li\\
  School of Mathematics and Statistics\\
  Shandong University of Technology\\
  Zibo, Shandong 255049, China\\
  \texttt{zyli@sdut.edu.cn}\\
  \And
  Masahiro Yamamoto\\
  Graduate School of Mathematical Sciences\\
  The University of Tokyo\\
  3-8-1 Komaba, Meguro-ku, Tokyo 153-8914, Japan\\
  \texttt{myama@ms.u-tokyo.ac.jp}\\
}

\begin{document}
\maketitle

\begin{abstract}
When considering fractional diffusion equation as model equation in analyzing anomalous diffusion processes, 
some important parameters in the model, for example, the orders of the fractional derivative or the source term, are often unknown, which requires one to discuss inverse problems to identify these physical quantities from some additional information that can be observed or measured practically. This chapter investigates several kinds of inverse coefficient problems for the fractional diffusion equation.
\end{abstract}

\keywords{Inverse coefficient problem\and time fractional partial differential equation\and uniqueness\and stability\and lateral boundary data\and Carleman estimate\and Dirichlet-to-Neumann map}
\MRsubject{35R30\and 35R11\and 65M32}

%%%%%%%%%%%%%%%%%%%%%%%%%%%%%%%%%%%%%%%%%%%%%%%%%%%%%%%%%%%%%%%%%%%%%
\section{Introduction}
\label{sec-intro-multi}

In this chapter, let $\OOO \subset \R^d$ be a bounded domain with 
smooth boundary 
$\ppp\OOO$ and let $0 < \alpha < 1$, if we do not specify.  
We consider fractional diffusion equations:
\begin{equation}
\label{equ-potential}
\pppa y(x,t) = \Delta y(x,t) + p(x)y(x,t), \quad 
x\in \OOO, \thinspace 0 < t < T,
\end{equation}
and
\begin{equation}
\label{equ-diffu_coef}
\pppa y(x,t) = \mbox{div}\, (p(x)\nabla y(x,t)), \quad
x\in \OOO, \thinspace 0 < t < T.
\end{equation}
Here and henceforth $\pppa v$ denotes the Caputo derivative given by 
$$
\pa_t^{\alpha}v(t) = \frac{1}{\Gamma(1-\alpha)}\int^t_0 (t-s)^{-\alpha}
\frac{dv}{ds}(s) ds.
$$ 
To each of \eqref{equ-potential} and \eqref{equ-diffu_coef}, we attach boundary and initial conditions:
\begin{equation}
\label{condi-boundary}
y(x,t) = h(x,t), \quad x\in \ppp\OOO, \thinspace 0 < t < T.
\end{equation}
and
\begin{equation}
\label{condi-initial}
y(x,0) = a(x), \quad x\in \OOO.
\end{equation}
We can consider other type of boundary conditions but we mainly attach the
Dirichlet boundary condition.

Given $p$, $\alpha\in (0,1)$, $h$ and $a$ in \eqref{equ-potential}-\eqref{equ-diffu_coef}, it is a forward problem that we find $y$ satisfying \eqref{equ-potential} or \eqref{equ-diffu_coef} with \eqref{condi-boundary} and \eqref{condi-initial} in suitable function spaces. Then \eqref{equ-potential}-\eqref{condi-boundary}-\eqref{condi-initial} and \eqref{equ-diffu_coef}-\eqref{condi-initial} are called initial-boundary value problems.
On the other hand,  for example, in \eqref{equ-diffu_coef}, $p(x)$ describes a spatial diffusivity coefficient 
in the media under consideration and so it is physically important how to choose $p(x)$.  When we make modelling, we have to determine $\alpha \in (0,1)$ and $p(x)$, $x \in \OOO$ in order that the solution $y$ to the initial-boundary value problem \eqref{equ-potential}-\eqref{condi-boundary}-\eqref{condi-initial} or \eqref{equ-diffu_coef}-\eqref{condi-initial} behaves appropriately: for instance, in the one-dimensional case $d=1$, we are requested to determine $p(x)$, $x \in \OOO$ such that $y(x_0,t) = \mu(t)$ for $0 < t < T$, where $\mu(t)$ is a prescribed target function or observation data function and an $x_0\in \OOO$ is a fixed monitoring point.  This is one inverse coefficient problem for a fractional diffusion equation which we survey in this chapter.

For more precise formulations of the inverse coefficient problems, restricted to simple equations such as \eqref{equ-potential} and \eqref{equ-diffu_coef}, we consider a solution $y(\alpha,p;h,a)=y(\alpha,p;h,a)(x,t)$ to
\begin{equation}
\label{equ-IBVP}
\left\{ 
\begin{alignedat}{2}
& \pppa y(x,t) = \Delta y(x,t) + p(x)y(x,t), &\quad& 
x\in \OOO, \thinspace 0 < t < T,      \\
& y(x,t) = h(x,t), &\quad& x\in \ppp\OOO, \thinspace 0 < t < T,\\
& y(x,0) = a(x), &\quad& x\in \OOO.   
\end{alignedat}
\right.           
\end{equation}
Here and henceforth we understand the solution $y(\alpha,p;h,a)$ in suitable classes (e.g., strong solutions), on which we can consult Gorenflo, Luchko and Yamamoto \cite{GLY15}, Kubica and Yamamoto \cite{KY}, Zacher \cite{Zacher} for example.

In this chapter, we mainly survey results on the determination of spatially varying coefficients $p(x)$. Let $\gamma \subset \ppp\OOO$ be a suitable subboundary and $\nu = \nu(x)$ be the outward unit normal vector to $\ppp\OOO$. We set $\ppp_{\nu}u = \nabla u\cdot \nu$.

In Sections 2 and 3, we can survey the following two major formulations of the inverse coefficient problems.

{\bf Section 2: 
A single or a finite time of observations of lateral boundary data.}
\\
Let $h$ and $a$ be appropriately chosen. Determine $p(x)$ and/or $\alpha \in (0,1)$ by 
$$
\ppp_{\nu}y(\alpha,p;h,a)(x,t), \quad x \in \gamma, \thinspace 
0 < t < T.
$$

{\bf Section 3: 
Infinitely many times of observations of lateral boundary data.}
\\
We fix $a=0$ in \eqref{equ-IBVP} for example.  We define the Dirichlet-to-Neuman map by 
$$
\Lambda_{\alpha,p}: h \mapsto \ppp_{\nu}y(\alpha,p;h,0)(x,t)
\vert_{\gamma\times (0,T)}
$$
with suitable domain space (e.g., $C^{\infty}_0(\gamma \times (0,T))$) and range space. Then determine $p(x)$ and/or $\alpha \in (0,1)$ by $\Lambda_{\alpha,p}$.
\\

Moreover we can consider the following
\\
{\bf Final overdetermination.} 
\\
We fix $h$, $a$.  Determine $p(x)$ and/or $\alpha$ by 
$$
y(\alpha,p;h,a)(x,T), \quad x\in \OOO.
$$
However, for the best of the authors' knowledge, there are no publications on the inverse coefficient problems with final observation data. 
\\

The mathematical subjects are {\bf Uniqueness} and {\bf Stability} and we survey them.
\\

The rest sections of this chapter are
\\
{\bf Section 4: Other related inverse problems.
\\
Section 5: Numerical works.}

In the final section, we are not limited to inverse coefficient problems, and we review numerical works on inverse problems, which include the lateral Cauchy problem and the backward problem in time, because many numerical articles are published rapidly also by practical demands and we think that even a tentative review may be useful for overlooking future researches.

We close this section with descriptions of available methodologies for inverse problems for fractional partial differential equations. As long as the cases of $\alpha=1$ and $\alpha=2$ are concerned, a comprehensive method for inverse problems with the formulation by a finite time of observations of lateral boundary data is based on Carleman estimates, which was created by Bukhgeim and Klibanov \cite{BuKl} and yields the uniqueness and the stability for inverse coefficient problems.  See also Beilina and Klibanov \cite{BK}, Bellassoued and Yamamoto \cite{BY}, Yamamoto \cite{Y2009} for example. Moreover for the formulation with infinitely many times of observations of lateral boundary data, an essential step is to construct special forms of solutions to differential equations under consideration, which are called complex geometric optics solutions, and Carleman estimates give flexible 
constructions (e.g., Bukhgeim \cite{Bu}, Imanuvilov and Yamamoto \cite{IY}).

We describe a Carleman estimate in a simplified form as follows. Let $Q \subset \R^{d+1}$ be a domain in the $(x,t)$-space, $\va = \va(x,t)$ a suitably chosen weight function and let $P$ be a partial differential operator. Then a Carleman estimate is stated by: there exist constants $C>0$ and $s_0 > 0$ such that 
$$
\int_Q s\vert u(x,t)\vert^2 e^{2s\va(x,t)} dxdt 
\le C\int_Q \vert Pu(x,t)\vert^2 e^{2s\va(x,t)} dxdt
$$
for all $s \ge s_0$ and $u \in C^{\infty}_0(Q)$.

Here we note that the Carleman estimate holds uniformly for all sufficient large $s>0$.  In other words, the constant $C>0$ is chosen independently of $s\ge s_0$ as well as $u$. The power of $s$ on the left-hand side can change according to partial differential equations.

For proving Carleman estimates, the indispensable step is the integration by parts:
\begin{equation}
\label{int_by_parts}
\ppp(yz) = y\ppp z + z\ppp y,
\end{equation}
where $\ppp = \frac{\ppp}{\ppp x_i}$ or $= \frac{\ppp}{\ppp t}$. Carleman estimates have been established for various equations such as the parabolic equation ($\alpha=1$) and the hyperbolic equation ($\alpha=2$) and we refer for example to Bellassoud and Yamamoto \cite{BY}, Isakov \cite{Is}, Yamamoto \cite{Y2009}. Here we do not intend any complete list of references on Carleman estimates.

However for fractional derivatives, convenient formulae such as \eqref{int_by_parts}, do not hold, so that we cannot prove Carleman estimates in general for fractional partial differential equations.

Thus for inverse problems for fractional cases, we have no comprehensive methodologies and there are not many  results on the uniqueness and the stability for the inverse problems in spite of the significance.

For the inverse coefficient problems, in Sections 2 and 3, we can mainly have the following strategies.
\begin{itemize}
\item
{\bf Representation of solutions to \eqref{equ-potential} by means of eigenfunctions.}\\
We extract information of the spectrum of the operator $\Delta + p(x)$ with the boundary condition. As for the representation of the solution, see, for example, Sakamoto and Yamamoto \cite{SY11}. In particular, in the one-dimensional case $d=1$, we can apply the Gel'fand-Levitan theory.
\item
{\bf Transforms in $t$.}\\
We can apply integral transoforms such as the Laplace transform, or the limits as $t \to \infty$.
\item
{\bf Reduction to partial differential equations with integer-orders.}\\
In special cases (e.g., $\alpha = \frac{1}{2}$), the reduction is possible, so that we can prove a Carleman estimate.
\end{itemize}

\section{A finite number of observations of lateral boundary data}

In this section, we present several results regarding inverse problems of determining the coefficients in the fractional diffusion equations from lateral Cauchy data.

\subsection{One-dimensioncal case by the Gel'fand-Levitan theory}

We start with a one-dimensional fractional diffusion equation:
\begin{equation}
\label{equ-CNYY}
\left\{ 
\begin{alignedat}{2}
&\pa_t^\al y(x,t) =\pa_x \left(D(x) \pa_x y(x,t)\right), &\quad&  0 <x < 1, \,
0<t<T, \\
&y(x,0)=\delta(x), &\quad&  0 < x < 1,\\
&\ppp_xy(0,t) = \ppp_xy(1,t) = 0, &\quad&  0 < t < T.
\end{alignedat}
\right.
\end{equation}
Here $\delta$ is the Dirac delta function.

We discuss the following
\begin{prob}
Determine the order $\al\in(0,1)$ of the time derivative and the diffusion coefficient $D(x)$ from boundary data $y(0,t)$, $0<t\le T$.
\end{prob}

We introduce the admissable set
$$
\cU:=\{ (\alpha,D);\ \al\in (0,1), D\in C^2[0,1], D>0 \mbox{ on } [0,1]\},
$$
and then have
\begin{thm}[Cheng, Nakagawa, Yamamoto and Yamazaki \cite{CNYY09}]
\label{thm-CNYY09}

Let $y(\alpha_1,D_1)$, $y(\alpha_2,D_2)$ be the weak solutions to \eqref{equ-CNYY} with respect to $(\alpha_1,D_1), (\alpha_2,D_2)\in\cU$, respectively. Then $y(\alpha_1,D_1)(0,t) = y(\alpha_2,D_2)(0,t)$, $0<t\le T$ with some 
$T>0$, implies $\alpha_1=\alpha_2$ and $D_1(x)=D_2(x)$, $0\le x\le 1$.
\end{thm}

By an argument similar to the above theorem, Li, Zhang, Jia and Yamamoto \cite{LZJY13} proved the uniqueness for determining the fractional order and the diffusion coefficient in the one-dimensional time-fractional diffusion 
equation with smooth initial functions by using boundary measurements. For $i=1,2$, let $\{ \va^{(i)}_n\}_{n\in \N}$ be the set of all the orthonormal eigenfunctions of $\frac{d}{dx}\left( D_i(x)\frac{d}{dx}\right)$ with the zero Neumann boundary condition $\frac{d\va}{dx}(0) = \frac{d\va}{dx}(1)=0$. 
 
\begin{thm}[Li, Zhang, Jia and Yamamoto \cite{LZJY13}]
\label{thm-LZJY13}
Suppose that $(\alpha_i,D_i) \in \cU$, and let $y(\alpha_i,D_i)$, $i=1,2$ satisfy
$$
\left\{ 
\begin{alignedat}{2}
&\pa_t^{\al_i} y(x,t)=\pa_x\left(D_i(x) \pa_x y(x,t) \right), &\quad&  
 0 < x < 1, \, 0 < t < T, \\
&y(x,0)=a(x), &\quad&  0 < x < 1,\\
&\ppp_xy (0,t) = \ppp_xy (1,t) = 0, &\quad&  0 < t < T.
\end{alignedat}
\right. 
$$
Suppose
$$
a\in H^4(0,1),\ a''(0)\neq0,\ a'(0)=a'(1)=0
$$
and
\begin{equation}
\label{condi-a_n}
\int^1_0 a(x)\varphi^{(1)}_n(x) dx \neq0,\quad \forall n\in\mathbb N \quad\mbox{ or } \quad \int^1_0 a(x)\varphi^{(2)}_n(x) dx \neq0,\quad \forall n\in\mathbb N.
\end{equation}
If $y(\alpha_1,D_1)(0,t) = y(\alpha_2,D_2)(0,t)$ and $y(\alpha_1,D_1)(1,t) = y(\alpha_2,D_2)(1,t)$, $0<t\le T$ with some $T>0$, then $\alpha_1=\alpha_2$ and $D_1(x)=D_2(x)$, $0 \le x \le 1$.
\end{thm}

Their result corresponds to classical results for the case $\alpha=1$ by Murayama \cite{Mu}, Murayama and Suzuki \cite{MuSu}.

As initial values, in Theorem 1, we have to choose exactly the Dirac delta function, while in Theorem 2, the initial value must satisfy very restricted non-degeneracy condition \eqref{condi-a_n}. One natural conjucture is that the uniqueness holds if $a \not\equiv 0$. It has not been proved even in the case of $\alpha=1$, However Pierce \cite{Pi} proved the uniqueness for the one-dimensional heat equation with non identically vanishing boundary data.
We can prove a similar result for the one-dimensional fractional diffusion equation.
\begin{thm}[Rundell and Yamamoto \cite{RY}]
\label{thm-RY}
Suppose that $0 < \alpha < 1$, $p_1, p_2 \in C[0,1]$, $p_1, p_2 < 0$ on $[0,1]$. Let $y(p_k)$ be sufficiently smooth and satisfy
$$
\left\{ 
\begin{alignedat}{2}
&\pa_t^{\al} y(x,t) =\pa_x^2 y(x,t) + p_k(x)y(x,t), &\quad & 0 < x < 1, \,
0 < t < T, \\
&y(x,0)=0, &\quad & 0 < x < 1,\\
&\ppp_xy (0,t) = 0, \quad \ppp_xy (1,t) = h(t), &\quad & 0 < t < T.
\end{alignedat}
\right.               
$$
If $h \not\equiv 0$ in $(0,T)$, then $y(p_1)(1,t) = y(p_2)(1,t)$, $0 < t < T$ yields  $p_1(x)=p_2(x)$, $0\le x \le 1$.
\end{thm}

Once that we establish the representation formula for the solution $y(p_k)$, the proof is similar to \cite{Pi}.  We can similarly prove also the uniqueness for $\alpha$ as well as $p(x)$ by the same data, but we omit the details.

For the one-dimensional time-fractional diffusion equation with spatially dependent source term:
\begin{equation}
\label{equ-JR}
\left\{ 
\begin{alignedat}{2}
& \partial_t^\alpha y(x,t) - \ppp_x^2y(x,t) + p(x)y(x,t)=f(x), &\quad& 
0<x<1, \, t>0,\\
&y(0,t)=y(1,t)=0, &\quad& t>0, \\
& y(x,0)=0, &\quad&  0<x<1.
\end{alignedat}
\right. 
\end{equation}
Jin and Rundell \cite{JR12} investigated other type of inverse coefficient problem, and gave the following uniqueness of the inverse problem.
\begin{thm}[Jin and Rundell \cite{JR12}]
\label{thm-JR12}
Let $y(p_i,f_j)(x,t)$ be the solution of \eqref{equ-JR} with $\{f_j\}_{j\in\mathbb N}$ and $p_i,i=1,2$ in $\{ p\in L^\infty(0,1); \,p\ge0,\|p\|\le M \}$. Here, $M>0$ is an arbitrarily fixed constant. We suppose that the set $\{f_j\}_{j\in\N}$ of input sources forms a complete basis in $L^2(0,1)$.
\begin{enumerate}
\item 
If further, $p_1=p_2$ on the interval $[1-\delta,1]$ for some $\delta\in(0,1)$ and $M<\delta\pi$, then there exists a time $t^*>0$ such that 
$$
\ppp_xy(p_1,f_j)(0,t) = \ppp_xy(p_2,f_j)(0,t),\quad j\in\N,
$$
for any $t >t^*$ implies $p_1=p_2$ on the interval $[0,1]$.
\item 
If $M<\pi$, then there exists a time $t^*>0$ such that
$$
\ppp_xy(p_1,f_j)(1,t) -  \ppp_xy(p_1,f_j)(0,t) = \ppp_xy(p_2.f_j)(1,t) -  \ppp_xy(p_2,f_j)(0,t),\quad j\in\N,
$$
for any $t>t^*$ implies $p_1=p_2$ on the interval $[0,1]$.
\end{enumerate}
\end{thm}
As for related results, see also Jin and Rundell \cite{JR15}.

\subsection{Integral transform}

For the general dimensional case, we refer to Miller and Yamamoto \cite{MY13} which applied an integral transform which is given in terms of the Wright function (e.g., Bazhlekova \cite{Ba}). For presenting the main result in \cite{MY13}, some notations and settings are needed. 

Let $y=y(\alpha,p)$ be a smooth solution to
$$
\left\{ 
\begin{alignedat}{2}
&\pa_t^{\al} y(x,t) = \Delta y(x,t) + p(x)y(x,t), &\quad & x\in\Omega, \,0 < t < T, \\
&y(x,0)=a(x), &\quad & x\in\Omega, \quad\mbox{if } 0<\alpha<1,\\
&y(x,0)=a(x), \pa_t y(x,0)=0, &\quad & x\in\Omega, \quad\mbox{if } 1<\alpha<2,\\
&y(x,t) = 0, &\quad & x\in\pa\Omega,\, 0 < t < T.
\end{alignedat}
\right.               
$$
Here, $\Omega\subset\mathbb R^d$, $d\ge1$ is bounded domain with smooth $\pa\Omega$. Let $\omega$ be a subdomain of $\Omega$ such that $\partial\omega\supset\partial\Omega$. We set
$$
\mathcal U_M
:=\{ p\in W^{1,\infty}(\Omega); \,p\le0 \mbox{ in } \Omega, \|p\|_{W^{1,\infty}(\Omega)} \le M, p|_{\omega}=\eta \}
$$
with a constant $M>0$ and a smooth function $\eta$, which are arbitrarily chosen. Moreover, the initial value $a$ of the diffusion system satisfies
\begin{equation}
\label{condi-delta_a}
a\in H^3(\Omega)\cap H_0^2(\Omega),\quad \Delta a\in H_0^1(\Omega)
\end{equation}
and
\begin{equation}
\label{condi-a>0}
a(x)>0,\quad x\in \overline{\Omega\setminus\omega}.
\end{equation}
Now we are ready to state the unique result from \cite{MY13}.
\begin{thm}[Miller and Yamamoto \cite{MY13}]
\label{thm-MY13}
Let $\alpha_1,\alpha_2\in(0,1)\cap(1,2)$, and the initial value satisfy \eqref{condi-delta_a} and \eqref{condi-a>0}, and let $p_1,p_2\in\mathcal U_M$:
\begin{enumerate}
\item 
If $y(\alpha_1,p_1)=y(\alpha_1,p_2)$ in $\omega\times(0,T)$, then $p_1=p_2$ in $\Omega$.
\item 
Moreover, let
$$
a\le0\mbox{ or } \ge0,\not\equiv0.
$$
\end{enumerate}
If $y(\alpha_1,p_1)=y(\alpha_2,p_2)$ in $\omega\times(0,T)$, then $\alpha_1=\alpha_2$ 
and $p_1=p_2$ in $\Omega$.
\end{thm}

Here the natural numbers for $\alpha_1,\alpha_2$ are excluded. However, using the existing uniqueness result for the parabolic inverse problem (e.g., pp 594--595 in \cite{Kli92}), Theorem 5 holds for all $0 < \alpha_1,\alpha_2< 2$.

As other proof of Theorem 5, we can argue as follows. By the formula on the Laplace transform of the Caputo 
derivative,
$$
L(\ppp_t^{\alpha}v)(s) := \int^{\infty}_0 e^{-s t} \pppa v(t) dt = s^\alpha L(v)(s) - s^{\alpha-1}v(0)
$$
for smooth $v$, we can apply also the Laplace transform of the fractional equations under consideration. We omit the details.

\subsection{Carleman estimates in restricted cases}

In Section 1,we have described no general Carleman estimates for fractional partial differential equations, unlike the classical partial differential equations. Thus the stability results for the inverse problem for the fractional 
diffusion equation with general fractional order from the lateral Cauchy data are rather limited. However, for some special orders such as $\alpha=1/2$, some Carleman estimates were proved and/or applied to inverse problems in Ren and Xu \cite{RX}, Xu, Cheng and Yamamoto \cite{XCY11}, and Yamamoto and Zhang \cite{YZ12}.

Yamamoto and Zhang \cite{YZ12} considered the following fractional diffusion equation with an initial condition and Cauchy data in a spatially one-dimensional case:
\begin{equation}
\label{equ-YZ}
\left\{ 
\begin{alignedat}{2}
& (\partial_t^{\frac12} -\partial_x^2)y(x,t) = p(x)y(x,t), &\quad&  0<x<1, \, 0<t<T,\\
& y(x,0)=a(x), &\quad& 0 < x < 1,\\
& y(0,t)=h_0(t),\quad  \partial_xy(0,t)=h_1(t), &\quad& 0< t < T,
\end{alignedat}
\right.
\end{equation}
where $a(x)$, $h_0(t), h_1(t)$ are given, and $a\not\equiv 0$ in $(0,1) \times(0,T)$ or $h_0, h_1 \not\equiv0$ in $(0,T)$, and established a Carleman estimate and proved the conditional stability in determining a coefficient $p(x)$ for a fractional diffusion equation. 

Following \cite{YZ12}, we let $t_0 \in (0,T)$ be arbitrarily fixed. We choose $x_0 >1$ such that $(x_0 - 1)^2 < \frac13$, and set
$$
d(x)=|x-x_0|^2,\ \psi(x,t)=d(x) - \beta (t-t_0)^2>0,\quad 0\le x \le 1, \, 0\le t \le T,
$$
where $\beta>0$ satisfies
$$
\sqrt{\frac{x_0^2}{\beta}} < \min\{ t_0,T-t_0 \}.
$$
For $x_0^2 > \varepsilon>(x_0-1)^2$, we set
$$
Q_\varepsilon := \{ (x,t)\in(0,1)\times(0,T); \psi(x,t)>\varepsilon \},
\quad \Omega_\varepsilon := Q_{\varepsilon}\cap \{t=t_0\}.
$$
Then we have the following.
\begin{thm}[Yamamoto and Zhang \cite{YZ12}]
For $k=1,2$, let $y(p_k)$ satisfy \eqref{equ-YZ} with $p_k$. Assume that there exist $M>0$ and sufficiently small $\varepsilon_0>0$ such that
\begin{align*}
&\|y(p_k)\|_{C^2([\varepsilon_0,T-\varepsilon_0];W^{2,\infty}(\Omega)\cap 
H^4(\Omega)) \cap C^3([\varepsilon_0,T-\varepsilon_0];L^2(\Omega))} 
+ \|y(p_k)\|_{C([0,T];L^\infty(\Omega))} \le M,\\
&\partial_x^jp_1(0) = \ppp_x^jp_2(0),\ j=0,1,\quad 
\|p_k\|_{W^{2,\infty}(\Omega)} \le M, \quad k=1,2.
\end{align*}
Moreover suppose
$$
y(p_k)(x,t_0) \ne0,\quad 0\le x\le 1, \quad \mbox{$k=1$ or $k=2$}.
$$
Then, for any $\varepsilon\in((x_0-1)^2, \frac13 x_0^2)$, there exists constants $C>0$ and $\theta \in(0,1)$ depending on $M$ and $\varepsilon$ such that
$$
\|p_1-p_2\|_{H^2(\Omega_{3\varepsilon}) } \le C\|(y(p_1) - y(p_2))(\cdot,t_0)\|_{H^4(\Omega_\varepsilon)}^\theta.
$$
\end{thm}

Unlike inverse coefficient problems for equations with natural number orders, one need two spatial data $y(\cdot,0)$ and $y(\cdot,t_0)$ to construct the unknown coefficients. By $(x_0-1)^2 < \ep < \frac{1}{3}x_0^2$, we note that 
$x_0 - \sqrt{3\ep} < 1$, that is, $(0, x_0 - \sqrt{3\ep}) \subset \subset (0,1)$.  Thus, in \cite{YZ12}, only in a subinterval $0 < x < x_0 - \sqrt{3\ep}$ of $[0,1]$, we can estimate $p_1-p_2$, which means that Theorem 6 gives the stability locally in $(0,1)$

Next, let $0 < t_0 < T$.  For a system
\begin{equation}
\label{equ-Ka}
\left\{ 
\begin{alignedat}{2}
& (\partial_t^{\frac12} -\partial_x^2)u(x,t) = g(x)\rho(x,t), &\quad& 
0<x<1, \, 0<t<T,\\
& u(x,0)=0, &\quad& 0 < x < 1,\\
& u(0,t)=u(1,t) = 0, &\quad&  0< t < T.
\end{alignedat}
\right.
\end{equation}
Kawamoto \cite{Ka1} discussed an inverse source problem: given $\rho$, determine $g=g(x)$, $0 < x < 1$ by data
$\ppp_x u(0,t)$, $\ppp_x u(1,t)$ for $0 < t < T$ and $u(x,t_0)$ for $0 < x < 1$ and estimated $g$ over the whole
interval $0 < x < 1$. We note that the system \eqref{equ-Ka} is a linearization of \eqref{equ-YZ} around $p=0$.  

For \eqref{equ-YZ}, his result easily yields the estimate $\Vert p_1 - p_2\Vert_{H^2(0,1)}$ for the inverse coefficient problem, which is in contrast with the local estimate of Theorem 6. We notice that \cite{Ka1} considered a more general elliptic operator but we omit the details.

It should be mentioned here that in general the regularity of the solution at $t = 0$ cannot be improved more than the assumptions in Theorem 6. For concentrating on the inverse problem, the discussion for the regularity of the solution to the initial-boundary value problem \eqref{equ-YZ} is omitted.

The proof of the Carleman estimate for the one-dimensional case with order $\alpha = \frac{1}{2}$ is done
by twice applying the Caputo derivative to convert the original fractional diffusion equation to a usual partial differential operator: $\partial_t - \partial_x^4$.Such arguments for general rational $\alpha$ is direct but 
extremely complicated, and so far it turns out not to work.

For fractional partial differential equations in general dimensions, a different method brings Carleman estimates:
Cheng, Lin and Nakamura \cite{CLN} for $\alpha=\frac{1}{2}$, Lin and Nakamura \cite{LiNa1} for $0 < \alpha < 1$ and
Lin and Nakamura \cite{LiNa2} for multi-term time fractional derivatives of orders $\in (0,1) \cap (1,2)$.  The Carleman estimate in \cite{LiNa1} produces the following uniqueness: if $y$ with some regularity satisfies
$$
\left\{ 
\begin{alignedat}{2}
&\pppa y = \Delta y(x,t) + \sum_{j=1}^d b_j(x,t)\ppp_jy + c(x,t)y, &\quad& x \in \OOO, \, 0 < t < T, \\
& y(x,0) = 0, &\quad& x\in \OOO, 
\end{alignedat}
\right.
$$
and
$$
y=0 \quad \mbox{in $\omega \times (0,T)$}
$$
with some subdomain $\omega \subset \OOO$, then 
$$
y=0 \quad \mbox{in $\OOO \times (0,T)$}.
$$
Here, we assume $b_j, c \in L^{\infty}(\OOO \times (0,T))$.

In general, we need to assume that $y(\cdot,0) = 0$ in $\OOO$, which is not requested for the unique continuation for the classical partial differenetial equations. 

Li and Yamamoto \cite{LiYa} proved that if $y \in L^{\infty}(0,T;H^2(\OOO))$ satisfies 
$$
\left\{ 
\begin{alignedat}{2}
&\pppa y(x,t) = \ppp_x^2y(x,t), &\quad& 0<x<1, \, 0 < t < T, \\
& y(x,t)=0, &\quad& x\in I, \, 0<t<T
\end{alignedat}
\right.
$$
with some non-empty open interval $I \subset (0,1)$, then 
$$
y=0 \quad \mbox{in $(0,1) \times (0,T)$}. 
$$

The applications of the Carleman estimates in \cite{CLN}, \cite{LiNa1} and \cite{LiNa2} require a transform of $(x,t)$ called the Holmgren transform.  Therefore for the inverse coefficient problem, the Holmgren transform does not keep the structure of the equation in \eqref{equ-IBVP}.  Moreover for the uniqueness, the zero initial condition is needed and so their Carleman estimates are not applicable to inverse source or inverse coefficient problems as we formulate.

We further review recent works on Carleman estimates. Kawamoto and Machida \cite{KaMa} considers a fractional transport equation:
\begin{align}
&\ppp_t^{\frac{1}{2}}y(x,v,t) + v\ppp_xy(x,v,t) + \sigma_t(x,v)y(x,v,t)
\label{equ-KM}\\
= &\sigma_s(x,v)   \int_V r(x,v,v')y(x,v',t) dv', \quad 
x \in \OOO := (0,\ell), \, v \in V:= \{ v_0 \le \vert v\vert \le v_1\},\nonumber
\end{align}
with
$$
y(x,v,0) = a(x,v), \quad x \in \OOO, \, v\in V
$$
and
$$
y(x,v,t) = h(x,v,t), \quad (x,v) \in \Gamma_+, \, 0 < t < T.
$$
Here, $v_0, v_1, \ell > 0$ are constants and we set 
$$
\Gamma_{\pm} := \{ (x,v) \in \ppp\OOO \times V; \,
\pm v < 0\, \mbox{at $x=0$}, \quad 
\pm v > 0\, \mbox{at $x=\ell$}\},
$$
where the double-signs correspond.  We note that $y(x,v,t)$ describes the density of some particles at the point
$x \in \OOO$ and the time $t$ with the velocity $v$. As for some physical backgrounds (see e.g., Machida \cite{Ma}), and the equation where $\ppp_t^{\frac{1}{2}}$ is replaced by $\ppp_t$ is called a radiative transport equation (e.g., Duderstadt and Martin \cite{Duder}).

By an idea similar to Xu, Cheng and Yamamoto \cite{XCY11}, reducing \eqref{equ-KM} to $\ppp_ty - v^2\ppp_x^2v$ with integral term, the article \cite{KaMa} established a Carleman estimate and proved the Lipschitz stability in determining $\sigma_t(x,v)$ and/or $\sigma_s(x,v)$ by $y(\cdot,\cdot,t_0)$ in $\Omega\times V$, $\pa_xy(\cdot,\cdot, 0)$ in $\OOO \times J$ and $y\vert_{\Gamma_+\times J}$, where $t_0\in (0,T)$ is arbitrarily fixed and $J \subset (0,T)$ is an arbitrary open subinterval including $t_0$.

Again let $\OOO \subset \R^d$ be a smooth bounded domain. Now we review Carleman estimates for a mixed type of fractional derivatives:
$$
\ppp_ty + \rho\pppa y
= \sum_{i,j=1}^d \ppp_i(a_{ij}(x)\ppp_jy)
+ \sum_{j=1}^d b_j(x,t)\ppp_jy + c(x,t)y.
$$
Here, $\rho>0$ is a constant and we assume $a_{ij} = a_{ji} \in C^1(\ooo{\OOO})$, $b_j, c \in L^{\infty}(\OOO\times (0,T))$ and the uniform ellipticity.

Kawamoto \cite{Ka2} proved a Carleman estimate for $\alpha = \frac{1}{2}$, and Huang, Li and Yamamoto \cite{HLY18} established Carleman estimates for both cases
\begin{enumerate}
\item
$0 < \alpha < \frac{1}{2}$.
\item
$\alpha$ is a rational number such that $0 < \alpha \le \frac{3}{4}$.
\end{enumerate}

See also Kawamoto and Machida \cite{KaMa18}.

We note that the Carleman estimate in \cite{Ka2} holds also for a more general equation
\begin{align*}
&\ppp_ty + \sum_{k=1}^{\ell} q_k(x,t)\ppp_t^{\alpha_k}y
= \sum_{i,j=1}^d \ppp_i(a_{ij}(x,t)\ppp_jy)\\
+ &\sum_{j=1}^d b_j(x,t)\ppp_jy + c(x,t)y, \quad 
0<\alpha_{\ell}< .... < \alpha_1 \le \frac{1}{2},
\end{align*}
with suitable conditions on $a_{ij}, b_j, c$ and see also \cite{HLY18} for possible generalizations.

According to the method by Bukhgeim and Klibanov \cite{BuKl}, their Carleman estimates yield the stability and the uniqueness for the inverse source problems under some suitable formulations, but we omit further details.

\section{Infinite many times of observations of lateral boundary data}

We fix $\la \in C^{\infty}[0,\infty)$ such that 
$$
\la(0) = \frac{d\la}{dt}(0) = 0
$$
and there exist constants $C>0$ and $\theta \in \left( 0, \frac{\pi}{2}\right)$ such that $\la$ can be analytically extended to $\{ z\in \C\setminus \{ 0\};\, \vert \mbox{arg} z\vert < \theta\}$ and 
$$
\left\vert \frac{d^2\la}{dt^2}(t)\right\vert = O(e^{Ct})
$$
as $t \to \infty$.  First, letting boundary values take the form of separation of variables, we consider an inverse coefficient problem with infinitely many times of observations .
 
We discuss the following initial-boundary value problem
\begin{equation}
\label{equ-LIY}
\left\{ 
\begin{alignedat}{2}
&\sum_{k=1}^{\ell} p_k(x)\pa^{\al_k}_t y(x,t) =\Delta y(x,t)+p(x)y(x,t),
&\quad& x\in \Om, \, 0<t<T,\\
&y(x,0)=0, &\quad& x \in \Om,\\
&y(x,t)= \la(t)h(x), &\quad& x\in \pa\Om, \, 0<t<T.
\end{alignedat}
\right. 
\end{equation}

We recall that $\nu$ is the outward unit normal vector to $\pa\Om$ and we denote $\ppp_{\nu}v = \na v\cdot\nu$. For $\ell \in \N$, we set $\bm\al= (\al_1,\cdots,\al_\ell) \in (0,1)^\ell$ where $\al_\ell < \al_{\ell-1} < \cdots < \al_1$.

We study
\begin{prob}
For $h\in H^{\f 3 2}(\pa\Om)$, we define the Dirichlet-to-Neumann map by
$$
\La(\ell,\bm \al,p_k,p)h:= \ppp_{\nu}y\Big|_{\pa\Sg}\in L^2(0,T;H^{\f 1 2}(\pa\Om)).
$$
Then we discuss whether $(\ell,\bm\al,p_k,p)$ is uniquely determined by the Dirichlet-to-Neumann map $\La(\ell,\bm\al,p_k,p):H^{\f32}(\pa\Om) \longrightarrow L^2(0,T;H^{\f12}(\pa\Om))$.
\end{prob}

As an admissible set of unknowns including numbers and coefficients, we set
\begin{align*}
\cU= & \{(\ell, \bm \al, p_k, p)\in\N\times(0,1)^\ell\times 
C^\infty(\ov\Om)^{\ell+1}, \, p_1 > 0, \\
& p_k \ge 0, \not\equiv 0, \, \mbox{$2\le k\le\ell$, $p\le0$ on $\ov\Om$} \}.
\end{align*}

The analyicity in $t$ of the solution to \eqref{equ-LIY} reduces the Dirichlet-to-Neumann map for \eqref{equ-LIY} by letting $t\to \infty$ to an inverse coefficient problem by elliptic Dirichlet-to-Neumann map. Therefore one can apply the existing uniqueness results.  

\begin{thm}[Li, Imanuvilov and Yamamoto \cite{LIY16}]
\label{thm-ip-DN_map}
Let $(\ell,\bm\al,p_k,p), (m,\bm \be,q_k,q) \in \cU$. Then 
$$
\La(\ell,\bm\al,p_k,p)h=\La(m,\bm\be,q_k,q)h,\quad\forall h
\in H^{\f32}(\pa\Om)
$$
implies that
$$
\ell=m,\quad \bm\al=\bm\be, \quad p_k=q_k, \quad 1\le k\le \ell, \quad
p=q \quad \mbox{in $\OOO$}.
$$
\end{thm}

In $\cU$, the regularity of $p, p_1, ...., p_\ell$ can be relaxed but we do not discuss here.  In particular, the result by Dirichlet-to-Neumann map in the two-dimensional case $d=2$ (e.g., Imanuvilov and Yamamoto \cite{IY}) yields a sharp uniqueness result where Dirichlet inputs and Neumann outputs can be restricted on an arbitrary subboundary and the required regularity of unknown coefficients is relaxed.

\begin{coro}[Li, Imanuvilov and Yamamoto \cite{LIY16}]
\label{buharik}
Let $\Om\subset \R^2$ be a bounded domain with smooth boundary $\pa\Om$ and $\gamma \subset\pa\Om$ be an arbitrarily given subboundary and let $r>2$ be arbitrarily fixed. Assume 
\begin{align*}
\mathcal{\wh U}
= \Big\{(\ell, \bm\al, p_k, p) \in \N\times(0,1)^{\ell}\times
(W^{2,\infty}(\Om))^{\ell}\times W^{1,r}(\Om);\\
\mbox{ $p_1 > 0$, $p_k \ge 0, \not\equiv 0$, 
$2\le k\le\ell$, $p\le0$ on 
$\ov\Om$} \Big\}.
\end{align*}
If $(\ell,\bm\al,p_k,p)$ and $(m,\bm\be,q_k,q) \in \mathcal{\wh U}$ satisfy
$$
\La(\ell,\bm\al,p_k,p)h = \La(m,\bm \be,q_k,q)h \quad \mbox{ on $\gamma$}
$$
for all $h \in H^{\f 3 2}(\pa\Om)$ with {\rm supp}\! $g \subset \gamma$, then $\ell=m$, $\bm\al=\bm\be$, $p_k=q_k$, $1\le k\le \ell$ and $p=q$ in $\OOO$.
\end{coro}

Next we review an inverse coefficient problem by Dirichlet-to-Neumann map at one shot.  Assuming that $0 < \alpha <1$ or $1<\alpha<2$, we discuss
\begin{equation}
\label{equ-DNmap}
\left\{ 
\begin{alignedat}{2}
& \rho(x)\pppa y(x,t) = \mbox{div}\, (a(x)\nabla y(x,t))
+ V(x)y(x,t), &\quad& x\in \OOO, \, 0 < t < T, \\
& y(x,t) = h(x,t), &\quad& x \in\ppp\OOO, \, 0<t<T, \\
& y(x,0) = 0, &\quad& x \in \OOO \quad \mbox{if $0<\alpha<1$}, \\
& y(x,0) = \ppp_ty(x,0) =  0, &\quad& x \in \OOO \quad \mbox{if $1<\alpha<2$}. 
\end{alignedat}
\right. 
\end{equation}
Let $0 < T_0 < T$ be arbitrarily fixed.  We assume
$$
\rho_k, a_k \in C^2(\ooo{\OOO}), \, > 0 \quad \mbox{on $\ooo{\OOO}$}, \quad
V_k \in L^{\infty}(\OOO)
$$
with $k=1,2$. We define the Dirichlet-to-Neumann map $\Lambda_{\rho,a,V,T_0}$ at one shot. Let $S_{\rm in}, S_{\rm out} \subset \ppp\OOO$ be relatively open non-empty subboundaries and $S_{\rm in} \cap S_{\rm out} \ne \emptyset$, $\ooo{S_{\rm in} \cup S_{\rm out}} = \ppp\OOO$.  We define $\Lambda_{\rho,a,V,T_0}$ by
$$
\Lambda_{\rho,a,V,T_0}h = a(\cdot)\ppp_{\nu}y(\cdot,T_0)\vert
_{S_{\rm out}},
$$
where $y$ is the solution to \eqref{equ-DNmap}. Then we have the following.
\begin{thm}[Kian, Oksanen, Soccorsi and Yamamoto \cite{KOSY}]
\label{KOSYY}
We assume that either of the following conditions is fulfilled.
\begin{enumerate}
\item
$\rho_1 = \rho_2$ in $\OOO$ and $\nabla a_1 = \nabla a_2$ on $\ppp\OOO$.
\item
$a_1 = a_2$ in $\OOO$ and there exists a constant $C_0>0$ such that $\vert \rho_1(x) - \rho_2(x)\vert \le C_0|{\rm dist}\, (x, \ppp\OOO)|^2$ for $x \in \OOO$.
\item
$V_1 = V_2$ in $\OOO$, $\nabla a_1 = \nabla a_2$ on $\ppp\OOO$, and there exists a constant $C_0>0$ such that 
$\vert \rho_1(x) - \rho_2(x)\vert \le C_0|{\rm dist}\, (x, \ppp\OOO)|^2$ for $x \in \OOO$.
\end{enumerate}
Then 
$$
\La_{a_1,\rho_1,V_1,T_0}h = \La_{a_2,\rho_2,V_2,T_0}h \quad 
\mbox{on $S_{\rm out}$}
$$
for all $h\in C_0(S_{\rm in} \times (0,T_0))\cap C^2([0,T]; H^{\frac32}(\pa\Omega))$ implies $(a_1,\rho_1,V_1) = (a_2,\rho_2,V_2)$ in $\OOO$.
\end{thm}

Thus the Dirichlet-to-Neumann map can identify at most two coefficients in the elliptic part under some additional information on $\nabla a_k$ on $\ppp\OOO$ or $\rho_1-\rho_2$ in $\OOO$. See Canuto and Kavian \cite{CK} for the case of $\alpha=1$. The proof of the theorem is similar to \cite{CK} and one essential step is the eigenfunction expansion of $y$ and so it is very difficult to discuss the inverse coefficient problem of determining also $A(x)$
for a non-symmtric equation $\rho\pppa y = \mbox{div}\, (a\nabla y) + A(x)\cdot \nabla y + V(x)y$.
\\

In view of degrees of freedom, we compare Theorems 7 and 8.
\\
{\bf Theorem 7}.
\begin{itemize}
\item
Degree of freedom of unknowns: $d$\\
More precisely, each unknown function depends on $d$ spatial variables $x_1, ..., x_d$.
\item
Degree of freedom of data set: $(d-1) + d = 2d-1$\\
More precisely, inputs are a set of $h$'s, each of which depends on $(d-1)$-spatial variables,
and the output is Neumann data on $\ppp\OOO \times (0,T)$, that is, $(d-1)+1 = d$.  Therefore the degree of freedom is one of the product set of the sets of inputs and outputs sets, so that it is $2d-1$
\end{itemize}
Therefore our inverse problem is overdetermining for $d\ge 2$ because $d < 2d-1$.

{\bf Theorem 8}.
\begin{itemize}
\item
Degree of freedom of unknowns: $d$\\
\item
Degree of freedom of data: $d + (d-1) = 2d-1$.
\\
More precisely, input is Dirichlet data on $\ppp\OOO \times (0,T)$, that is,
the degree is $(d-1) + 1 = d$, while the output is the 
corresponding Neumann data at fixed time $T_0$, that is,
$d-1$. 
\end{itemize}
Therefore the formulation is equally overdetermining.
\\
\vspace{0.2cm}

So far we assume that the order $\alpha$ of the derivative is constant and is mainly in $(0,1)$. However, in some complex media, the presence of heterogeneous regions may cause variations of the permeability in different spatial positions, and in this case, the variable order time-fractional model is more relevant for describing the diffusion process (e.g., \cite{SCC}). 

Let us turn to considering an inverse problem of determining the variable order in $x$ and related coefficients.
We introduce several necessary notation and settings. We are given constants $0 < \alpha_0 < \alpha_M$, $0< \rho_0 < \rho_M$ and two functions $\alpha\in L^\infty(\Omega)$  and $\rho\in L^\infty(\Omega)$  satisyfing
\begin{equation}
\label{condi-alpha(x)}
0 < \alpha_0 \le \alpha(x) \le \alpha_M < 1\mbox{ and } 0 <\rho_0 \le \rho(x) 
\le \rho_M < \infty,\quad  x\in\Omega. 
\end{equation}
We consider an initial-boundary value problem for a space-dependent variable order fractional diffusion equation
\begin{equation}
\label{equ-SKY}
\left\{
\begin{alignedat}{2}
&\rho(x) \partial_t^{\alpha(x)} y(x,t) = \Delta y(x,t) +q(x)y(x,t), &\quad&
 x \in \Omega, \, t>0,\\
&y(x,0)=a(x), &\quad & x\in\Omega,\\
&y(x,t)=\la(t) h(x), &\quad & x \in \partial\Omega, \, t>0
\end{alignedat}
\right.
\end{equation}
with suitable $h$. Henceforth by $y_h$ we denote a unique solution in $C([0,\infty);H^2(\Omega))$ to \eqref{equ-SKY}. Given non-empty relatively open subboundaries $S_{\mathrm{in}}$ and $S_{\mathrm{out}}$ of $\partial\Omega$, we introduce the following boundary operator
$$
\Lambda(\alpha,\rho,q) : h \in \mathcal H_{\mathrm{in}} \mapsto 
\partial_\nu y_h(\cdot,t)|_{S_{\mathrm{out}}},\quad 0<t < \infty,
$$
where 
$$
\mathcal H_{\mathrm{in}}:=\{ h\in H^{3/2}(\partial\Omega); \,
\mathrm{supp}\ h \subset S_{\mathrm{in}} \}. 
$$

Kian, Soccorsi and Yamamoto \cite{KSY17} analysed the uniqueness in the inverse problem of determining simultaneously the variable order $\alpha$, coefficients $\rho$ and $q$ of the diffusion equation in \eqref{equ-SKY}
from the knowledge of the boundary operators $\{\Lambda(\alpha,\rho,q)(t_n);n\in\mathbb N\}$ associated with a time sequence $t_n$, $n\in\mathbb N$ fulfilling
\begin{equation}
\label{condi-tn}
\mbox{the set $\{t_n; \, n\in\mathbb N\}$ has an accumulation point in $(0,\infty)$. }  
\end{equation}

Moreover, $\Omega$, $S_{\mathrm{in}}$ and $S_{\mathrm{out}}$ are assumed to satisfy the following conditions.

{\bf Case I: } Spatial dimension $d=2$.\\
It is required that $\partial\Omega$ is composed of a finite number of smooth closed contours. In this case, choose $S_{\mathrm{in}}=S_{\mathrm{out}}=:\gamma$, where $\gamma$ is any arbitrary non-empty relatively open subset of $\partial\Omega$, and the set of admissible unknown coefficients reads
$$
\mbox{ $\mathcal U
:=\{ (\alpha,\rho,q)$; $\alpha,\rho\in W^{1,r}(\Omega)$ fulfill 
\eqref{condi-alpha(x)} and $q>0,\in W^{1,r}(\Omega)$\} },
$$
where $r\in(2,\infty)$.
\\
{\bf Case II: $d\ge3$}.\\
Let $x_0\in\mathbb R^d$ outside the convex hull of $\overline\Omega$, and assume that
$$
\{ x\in\partial\Omega; (x-x_0)\cdot \nu \ge0 \}\subset S_{\mathrm{in}} 
\mbox{ and } 
\{ x\in\partial\Omega; (x-x_0)\cdot \nu \le0 \}\subset S_{\mathrm{out}} 
$$
Furthermore we define the set of admissible unknown coefficients by
$$
\mbox{ $\mathcal V:=\{ (\alpha.\rho,q)$; $\alpha\in L^\infty(\Omega)$, 
$\rho\in L^\infty(\Omega)$ fulfill \eqref{condi-alpha(x)} and 
$q>0,\in L^\infty(\Omega)$\} }.
$$

The uniqueness results for the inverse coefficients problem are as follows.
\begin{thm}[Kian, Soccorsi and Yamamoto \cite{KSY17}]
Let $t_n$, $n\in\mathbb N$ fulfill \eqref{condi-tn}. We assume Case (I) or (II).  Let 
$$
(\alpha_i,\rho_i,q_i) \in \mathcal{U}, \quad i=1,2 \quad
\mbox{in Case I}
$$
and
$$
(\alpha_i,\rho_i,q_i) \in \mathcal{V}, \quad i=1,2 \quad
\mbox{in Case II.}
$$
If
$$
\Lambda(\alpha_1,\rho_1,q_1)(t_n)=\Lambda(\alpha_2,\rho_2,q_2)(t_n),\quad 
n\in\mathbb N, 
$$
then $(\alpha_1,\rho_1,q_1)=(\alpha_2,\rho_2,q_2)$.
\end{thm}
\section{Other related inverse problems}

\subsection{Determination of time varying coefficients}

For 
$$
\begin{cases}
\partial_t^\alpha y(x,t) - a(t)\ppp_x^2y(x,t) = 0, \quad & 0 < x < 1, \, 0<t<T, \\
y(x,0)=0, \quad & 0 < x < 1,\\
y(0,t)=h(t),\ y(1,t)=0, \quad & 0 < t < T,
\end{cases}
$$
where $\alpha\in(0,1)$ and $T>0$, $h$ are given, the article Zhang \cite{Z16} proved the uniqueness for an inverse problem of recovering $a(t)$ from additional boundary data
$$
-a(t) \ppp_xy(0,t),\quad 0<t<T.
$$

Moreover we can refer to Fujishiro and Kian \cite{FuKi}. In the case where unknown coefficients are dependent only on $t$, we can expect that the pointwise data of the solution at monitoring points over a time interval can guarantee the stability as well as the uniqueness.  Indeed we often reduce such an inverse problem to Volterra equation of the second kind.

Janno \cite{Ja} discussed an inverse problem of determining an order and a  time varying function in an integral term, while Wang and Wu \cite{WaWu} studied an inverse problem of determining two varying functions in an integral term and a source term.

\subsection{Determination of nonlinear terms}

Luchko, Rundell, Yamamoto and Zuo \cite{LRYZ} considers
$$
\left\{ 
\begin{alignedat}{2}
& \pppa y(x,t) = \Delta y(x,t) + f(y(x,t)), &\quad& x \in\OOO \subset 
\R^d, \, 0<t<T, \\
& \ppp_{\nu}y(x,t) = h(x,t), &\quad& x\in \ppp\OOO, \, 0 < t < T, \\
& y(x,0) = a_0, &\quad& x \in \OOO,
\end{alignedat}
\right.
$$
where $a_0$ is a constant.  The article \cite{LRYZ} discussed the determination of semilinear term $f$ by extra data $y(x_0,t)$, $0 < t < T$ at a fixed point $x_0 \in \ooo{\OOO}$ to establish a uniqueness result for the inverse problem within suitable admissible sets of $f$'s and provided a numerical method for reconstructing $f$.

Rundell, Xu and Zuo \cite{RXZ} studied a similar inverse problem of determining a nonlinear term $f$ in the boundary condition $-\ppp_xy(1,t) = f(y(1,t))$, $0 < t < T$ for the one-dimensional fractional diffusion equation. See Tatar and Ulusoy \cite{TU2017} for an optimization method for reconstructing $a$ in 
$$
\pppa y(x,t) = \mbox{div}\, (a(y(x,t))\nabla y(x,t)) + F(x,t), \quad 
x \in \OOO \subset \R^d, \, 0<t<T.
$$

Finally, we note that Janno and Kasemets \cite{JK2} considered an initial-boundary value problem for 
$$
\pppa y = \sum_{i,j=1}^d \ppp_i(a_{ij}(x)\ppp_jy) 
+ \sum_{j=1}^d b_j(x)\ppp_jy + p(y,x,t)g(x)
+ q(y,x,t)
$$
for $x \in \OOO$ and $0 < t < T$, and established the uniqueness in determining $g(x)$ by extra data 
$$
\int^T_0 y(x,t) \mu(t) dt
$$
with suitable weight $\mu(t)$.

Lopushanska and Rapita  \cite{LR} discussed an inverse problem of determining $r(t)$ in a semilinear fractional telegraph equation
$$
\pa_t^\alpha y+ r(t) \pa_t^\beta y = \Delta y + F(x,t,y,\pa_t^\beta y)
$$
where $1<\alpha<2$ and $0<\beta<1$ by means of data
$$
\int_\Omega y(x,t)w(x)dx
$$
with some weight function $w(x)$. See also Ismailov and \c{C}i\c{c}ek \cite{IC} for a similar inverse problem.

\section{Numeical works}
The literature on numerics increases up rapidly for inverse coefficient problems. Here, we refer to only a few papers: Bondarenko and Ivaschenko \cite{BS}, Li, Zhang, Jia and Yamamoto \cite{LZJY13}, Sun, Li and Jia \cite{SLJ}, Sun and Wei \cite{SW}. The publications so far are limited to individual inverse problems, and we expect more comprehensive numerical works.

Next, we briefly survey numerical works for other types of inverse problems.

\subsection{Lateral Cauchy problems}

In the lateral Cauchy problem, we are requested to determine $y(x,t)$ satisfying
\begin{equation}
\label{equ-cauchy}
\left\{ 
\begin{alignedat}{2}
& \pppa y(x,t) = \Delta y(x,t), &\quad& x \in \OOO\subset \R^d, 
\, 0<t<T,\\
&y(x,t) = h_0(x,t), \\
& \ppp_{\nu}y(x,t) = h_1(x,t), &\quad& x\in \Gamma \subset \ppp\OOO,
\, 0 < t < T
\end{alignedat}
\right.                 
\end{equation}
with given $h_0, h_1$.  Here, $\gamma \subset \ppp\OOO$ is a sub-boundary and the initial values are also unknown.  The lateral Cauchy problem is closely related to the unique continuation which is expected to be proved by Carleman estimates. Various numerical methods are available: Li, Xi and Xiong \cite{LXX}, Murio \cite{Mu2007, Mu2008, Mu2009}, Qian \cite{Qi}, Zheng and Wei \cite{ZW2010, ZW2012}.

We point our that theoretical researches for the uniqueness and the conditional stability in determining $y$ satisfying \eqref{equ-cauchy} are not satisfactorily made.

\subsection{Backward problems in time}

We consider 
\begin{equation}
\label{equ-back}
\left\{ 
\begin{alignedat}{2}
& \pppa y(x,t) = \sum_{i,j=1}^d \ppp_i(a_{ij}(x)\ppp_jy), &\quad& 
x\in \OOO, \, 0<t<T, \\
& y(x,t)=0, &\quad& x \in \ppp\OOO, \, 0<t<T.
\end{alignedat}
\right. 
\end{equation}
Here, $a_{ij} = a_{ji} \in C^1(\ooo{\OOO})$ and there exists a constant $\mu_0 > 0$ such that $\sum_{i,j=1}^d a_{ij}(x)\xi_i\xi_j \ge \mu_0\sum_{j=1}^d \xi_j^2$ for all $x \in \ooo{\OOO}$ and $\xi_1,..., \xi_d \in \R$.

We can discuss a more general ellitpic operator,  but it suffices to consider \eqref{equ-back}.

We consider the following.

{\bf Backward problem in time.}\\
Given $b\in L^2(\OOO)$, solve $y(x,t)$, $x \in \OOO$, $0<t<T$ satisfying \eqref{equ-back} and 
\begin{equation}
\label{obser-final}
y(x,T) = b(x), \quad x \in \OOO. 
\end{equation}

In the case of $\alpha=1$, it is well-known that this problem has no solutions in general, and is ill-posed. However, in the case of $0 < \alpha < 1$, the situation is completely different, and Sakamoto and Yamamoto \cite{SY11} proves that for given $b \in H^2(\OOO) \cap H^1_0(\OOO)$ there exists a unique solution $y \in L^2(0,T;H^2(\OOO)\cap H^1_0(\OOO)) \cap C([0,T]; L^2(\OOO))$ satisfying \eqref{equ-back} and \eqref{obser-final}, and there exist constants $C_1, C_2 > 0$ such that
$$
C_1\Vert b\Vert_{H^2(\OOO)} \le \Vert y(\cdot,0)\Vert_{L^2(\OOO)} \le C_2\Vert b\Vert_{H^2(\OOO)}
$$
for all $b \in H^2(\OOO)\cap H^1_0(\OOO)$.

That is, the backward problem in time for fractinal diffusion equations with order $\alpha \in (0,1)$ is well-posed in the sense of Hadamard if we strengthen the regularity of data $y(\cdot,T)$ in $H^2(\OOO) \cap H^1_0(\OOO)$, which means that the smoothing effect is with only $2$ by order of Sobolev spaces.  The smoothing in the case of $0 < \alpha < 1$ is much weaker than the case of $\alpha=1$. Thus we can expect more stable reconstruction of initial value by means of the final value $y(\cdot,T)$. In the case of $\alpha=1$, we can prove an inequality of Carleman type where the weight is given in the form $e^{\la t}$ with large constant $\la > 0$, and establish the conditional stability for more general parabolic equations.  However, for fractional partial differential equations, again by the lack of convenient formulae of integration by parts, we cannot prove such a Carleman estimate, so that the theoretical researches for the backward problem are very restricted.

We refer to numerical works:
\\
Liu and Yamamoto \cite{JJLY}, Tuan, Long and Tatar \cite{TLT}, Wang and Liu \cite{WL2012, WL2013}, Wang, Wei and Zhou \cite{WWZ}, Wei and Wang \cite{WeiZ}, Xiong, Wang and Li \cite{XWL2012}, Yang and Liu \cite{YL2013}.

%%%%%%%%%%%%%%%%%%%%%%%%%%%%%%%%%%%%%%%%%%%%%%%
\bibliographystyle{unsrt}

\end{document}